\documentclass[a4paper,12pt]{article}
\usepackage[english]{babel}
\usepackage{amsmath}
\usepackage{amssymb}
\usepackage{tabls}
\usepackage[pdftex]{graphicx}
\usepackage{calc}
\usepackage{xcolor}
\usepackage[T1]{fontenc}
\usepackage[latin1]{inputenc}
\title{Abelian Group Clifford Algebras}
\author{Tim Neijens\\University of Antwerp\\ \texttt{tim.neijens@gmail.com}}

\newcommand{\Z}{\mathbb{Z}}

\newcommand{\blok}{\hfill \Box}

\newtheorem{defi}{Definition}[section]

\newtheorem{prop}[defi]{Proposition}

\newtheorem{st}[defi]{Theorem}

\begin{document}
\maketitle
\begin{abstract}
In \cite{CVO} generalized Clifford algebras were introduced via Clifford representations; these correspond to projective representations of a finite group (Abelian), $G$ say, such that the corresponding twisted group ring has minimal center.  The latter then translates to the fact that the corresponding $2$-cocycle allows a minimal (none!) number of ray classes and this forces a decomposition of $G$ in cyclic components in a suitable way, cf. \cite{Z}.  In this small paper, I will provide a way to represent an Abelian Group Clifford Algebra using a matrix, and then give a way to calculate whether or not the center is trivial.
\end{abstract}

\section{Preliminaries and Notation}

{\defi A \textbf{projective representation of a finite group $G$} is a group morphism $T: G \rightarrow PGL_n(k)$ with $k$ some field or a commutative connected ring ($0$ and $1$ are the only idempotents).}\\

Such projective representations define algebra homomorphisms (by $k$-linearly extending $T$) $kG^\alpha \rightarrow M_n(k)$ where $kG^\alpha$ is the twisted group ring over $k$ with respect to some $2$-cocycle.  It is defined as
\[kG^\alpha = \bigoplus_{\sigma \in G} k u_\sigma,\]
and multiplication is defined by $u_\sigma u_\tau = \alpha(\sigma, \tau)u_{\sigma\tau}$ for $\sigma,\tau \in G$.

{\prop (see \cite{CVO}) Let $R$ be a commutative connected (only idempotents are $0$ and $1$), $G$ a (finite) group and $\alpha$ a $2$-cocycle, then the center $Z(RG^\alpha)$ is freely generated over $R$ by the ray class sums.  If $G$ is Abelian then $Z(RG^\alpha) = R(G_\textup{reg})^\alpha$.}\\
Let $G$ be a finite Abelian group.  We can decompose this group as follows:
\[G \cong G_{p_1} \oplus \ldots \oplus G_{p_t},\]
where
\[G_{p_v} \cong \left(\frac{\Z}{p_v^{n_1}}\right)^{m_1}\oplus \ldots \oplus \left(\frac{\Z}{p_v^{n_{k_v}}}\right)^{m_{k_v}},\]
with $n_1 > n_2 > \ldots > n_{k_v}$.\\
Suppose for now that $t=1$.  (It will become clear that there is no loss of generality by doing this.)  $G_p$ is generated by
\[
\begin{array}{*{20}c}
   {e_1 (n_1 ),} &  \ldots  & {,e_{m_1 } (n_1 )}  \\
    \vdots  & {} & {}  \\
   {e_1 (n_k ),} &  \cdots  & {,e_{m_k } (n_k )},  \\

 \end{array} 
\]
and denote $e_i(n_j)$ by $e_{ij}$, i.e. $e_{ij}$ has exact order $p^{n_j}$.\\
Now let $R$ be a commutative ring which contains a primitive $(p^{n_1})$-th root of unity and no idempotents other than $0$ and $1$.\\
Consider a $2$-cocycle $\phi: G \times G \rightarrow U(R) \in Z^2(G,U(R))$.  We can associate the map
\[f_\phi:G \times G \rightarrow U(R):(A,B)\mapsto\frac{\phi(A,B)}{\phi(B,A)}.\]
Since $G$ is Abelian, $f_\phi$ is a multiplicatively antisymmetric bipairing.  This implies that $f_\phi$ is completely determined by the images of $e_{ij}$.  Moreover, since the $e_{ij}$ have finite order, so do the images.  In other words:
\[f_\phi(e_{ij}, e_{rs})=\omega^x,\]
where $\omega$ is a $[p_{n_j}, p_{n_s}]$-th root of unity.  ($[a,b] = \gcd(a,b)$).  Remember that we assumed only one prime $p$.  For a projective representation $T$, $T(e_{ij})T(e_{rs})= f(e_{ij}, e_{rs})T(e_{rs})T(e_{ij})$.  This means that generators corresponding to different primes commute and therefore we get a tensor of all those components, each component corresponding to a different prime.  In other words, if $A_G$ is the Abelian group Clifford algebra corresponding to a group G, then we find for $G \cong G_{p_1} \oplus \ldots \oplus G_{p_t}$:
\[A_G \cong A_{G_{p_1}}\otimes \ldots \otimes A_{G_{p_t}}.\]
This means that, like said in the beginning of the paragraph, we only need to study the case that $G \cong G_p$.

\section{The Matrix Method}

Given a certain $2$-cocycle $\phi \in Z^2(G, U(R))$ we construct a matrix $A_\phi$ with entries $x_{ijrs}$ determined by
\[f(e_{ij}, e_{rs})= \omega_{[p^{n_j}, p^{n_s}]}^{x_{ijrs}}.\]
These entries are in the additively defined cyclic groups.  So we find:
\[
A_\phi=\left[ {\begin{array}{*{20}c}
   {A_{11} } & {A_{12} } &  \cdots  & {A_{1k} }  \\
   {A_{21} } & {A_{22} } &  \cdots  & {A_{2k} }  \\
    \vdots  & {} & {} &  \vdots   \\
   {A_{k1} } & {A_{k2} } &  \cdots  & {A_{kk} }  \\

 \end{array} } \right],
\]
where
\[A_{ij}\in M_{m_i \times m_j}\left(\frac{\Z}{p^{n_j}\Z}\right).\]
By construction, $A_\phi$ is antisymmetric.\\

If we want to minimize the center, we need to minimize the number of elements that commute with the generators of our projective representation $A$.  (We will from now on identify $t_{ij}$ with $T(e_{ij})$ for a projective representation $T$).  We will denote $f_\phi$ as $f$ if we fix $\phi$.  Let
\[b = t_{11}^{\alpha_{11}}\cdot \ldots \cdot t_{1m_1}^{\alpha_1 m_1}\cdot t_{21}^{\alpha_{21}}\cdot \ldots \cdot t_{k m_k}^{\alpha_{k m_k}}\in Z(A).\]
Then we find $\forall i,j$
\begin{eqnarray}
\ &\ & t_{ij}\cdot b \cdot t_{ij}^{-1} = b \label{cen12}\\
\ &\Rightarrow& f^{\alpha_{11}}(e_{ij}, e_{11})\cdot \ldots \cdot f^{\alpha_{k m_k}}(e_{ij}, e_{k m_k})\cdot b = b \\
\ &\Rightarrow& \omega_{[p^{n_j}, p^{n_1}]}^{\alpha_{11} x_{ij11}}\cdot \ldots\cdot \omega_{[p^{n_j}, p^{n_k}]}^{\alpha_{k m_k} x_{ijk m_k}}\cdot b = b. \label{cen13}
\end{eqnarray}
We want to minimize the number of solutions $(\alpha_{11}, \ldots, \alpha_{k m_k})$ to this equation.  We can express this minimality condition by using the matrix $A_\phi$ defined before.  We use this matrix to calculate the exponents of the $\omega$ appearing in (\ref{cen13}).  If we want the left factor in (\ref{cen13}) to be $1$, we have to calculate the product of the $\omega$ appearing in {\ref{cen13}}.  We do this by using a primitive $(p^{n_1})$-th root of unity, let us call it $\zeta$.  We then get as the exponent of $\zeta$:
\begin{equation}
p^{n_1-n_j}\cdot\alpha_{11} x_{ij11}+\ldots+p^{n_1-\textup{min}(n_j,n_k)}\alpha_{k m_k} x_{ijk m_k}. \label{cen14}
\end{equation}
$\zeta$ to that power must be equal to $1$ and that only happens if the exponent is $0 \bmod p^{n_1}$.  So we construct from $A_\phi$ a new matrix $\widetilde A_\phi$ by multiplying each block $A_{ij}$ with $p^{n_1-\textup{min}(n_i,n_j)}$.  The end result is that
\[\widetilde A_\phi \in M_{m_1 + \ldots +m_k}\left(\frac{\Z}{p^{n_1}\Z}\right).\]
Since the condition is now focused on the exponents, we are working with the additive groups and the standard ring structure on $\Z/n\Z$.  We will exploit this by using the matrix product with our newly formed matrix $\widetilde A_\phi$ to find solutions or the lack thereof.\\
Let $g = (g_1, \ldots, g_k)\in G$ where $g_i$ is an $m_i$-dimensional vector with entries chosen in representatives in $\Z/p^{n_1}\Z$ of elements in $\Z/p^{n_i}\Z$.  So each entry in $g$ modulo $p^{n_i}$ determines an element in the appropriate $\Z/p^{n_i}\Z$.  Now we can calculate the matrix product $\widetilde A_\phi \cdot g$.  If this product is equal to the $0$-vector, then we have found a solution $g$ to (\ref{cen14}).  In order to minimize the number of solutions of (\ref{cen14}), we have to maximize the 'rank' of $\widetilde A_\phi$.

\subsubsection{Example}

To make things clear, an actual example of how to solve the equations presented above.  We first put the matrix in standard form, meaning that we construct $\widetilde A_\phi$ from $A_\phi$.  The solutions $\widetilde X$ to $\widetilde A_\phi \widetilde X = 0$ are then reduced component per component with the corresponding prime power.  We can find those solutions using simple invertible row operations.  Let $G \cong (\Z_9)^2 \oplus (\Z_3)^2$:
\[
A_\phi   = \left( {\begin{array}{*{20}c}
   0 \hfill & 1 \hfill &\vline &  1 \hfill & 1 \hfill  \\
   8 \hfill & 0 \hfill &\vline &  2 \hfill & 2 \hfill  \\
\hline
   2 \hfill & 1 \hfill &\vline &  0 \hfill & 1 \hfill  \\
   2 \hfill & 1 \hfill &\vline &  2 \hfill & 0 \hfill  \\

 \end{array} } \right) \to \widetilde A_\phi   = \left( {\begin{array}{*{20}c}
   0 \hfill & 1 \hfill &\vline &  3 \hfill & 3 \hfill  \\
   8 \hfill & 0 \hfill &\vline &  6 \hfill & 6 \hfill  \\
\hline
   6 \hfill & 3 \hfill &\vline &  0 \hfill & 3 \hfill  \\
   6 \hfill & 3 \hfill &\vline &  6 \hfill & 0 \hfill  \\

 \end{array} } \right) \to \left( {\begin{array}{*{20}c}
   1 \hfill & 0 \hfill &\vline &  0 \hfill & 0 \hfill  \\
   0 \hfill & 1 \hfill &\vline &  0 \hfill & 0 \hfill  \\
\hline
   0 \hfill & 0 \hfill &\vline &  6 \hfill & 0 \hfill  \\
   0 \hfill & 0 \hfill &\vline &  0 \hfill & 3 \hfill  \\

 \end{array} } \right)
\]
and thus
\[
\left\{ {\begin{array}{*{20}c}
   {\tilde x_1  = 0} \hfill  \\
   {\tilde x_2  = 0} \hfill  \\
   {\tilde x_3  \in \left\{ {0,3,6} \right\}} \hfill  \\
   {\tilde x_4  \in \left\{ {0,3,6} \right\}} \hfill  \\

 \end{array} } \right. \Rightarrow \left\{ {\begin{array}{*{20}c}
   {x_1  = 0\bmod 9}  \\
   {x_2  = 0\bmod 9}  \\
   {x_3  = 0\bmod 3}  \\
   {x_4  = 0\bmod 3}  \\

 \end{array} } \right.
\]

\section{Diagonal Block Theorem}

In this section we prove a theorem providing the necessary conditions on $A_\phi$ and so on the commutation relation of a projective representation $T$ to form an Abelian group Clifford algebra.

{\st \label{cen15} $R(G)^\phi$ has trivial center (i.e. $Z(A)=R$) if and only if each block $A_{ii}$ in $A_\phi$ is invertible.}\\
\textbf{Proof} Since every matrix on the diagonal is invertible in its respective ring, we can use row operations to put those blocks into diagonal form.  We can 'create zeroes' using the entries in the now diagonalized blocks.  Fix such a diagonal block and as such a ring $\Z/p^s \Z$.  The entries in the blocks are invertible if we consider them in their 'original' ring (before we go from $A_\phi$ to $\widetilde A_\phi$) and the entries above or below these can be considered as in $\Z/p^s \Z$ by construction of the matrix $A_\phi$.  We can use row operations to eliminate the corresponding variables.  Repeating this we get a diagonal matrix.  This matrix has only the $0$-vector as solution, modulo the corresponding primes.\\
In order to establish the other direction, suppose that one or more matrices are not invertible in their respective rings.  Take the first such matrix, say $A_{ii}$.  Then consider $\widetilde A_\phi$ and do the row operations in $A_{ii}$ that express the linear dependence of the entries.  Let us say row $s$ is the first row in $\widetilde A_\phi$ that is linearly dependent of the other rows in $A_{ii}$.  Up to row $s$ we can put the matrix in diagonal form.  This means that either $x_s$ will be a fully random variable (the case when the rows below row $s$ have zeroes on column $s$), or we get an entry from a row from an $A_{jj}$ where $j>i$.  This entry is not invertible in the ring corresponding to $A_{ii}$ and therefore will not yield one unique solution in the corresponding ring.  So $A_{ii}$ must be invertible. $\blok$

\end{document}